\documentclass[12pt]{article}
  \usepackage{amsfonts,amsmath}
  \usepackage[square,sort&compress,comma,numbers]{natbib} 

  \usepackage{todonotes} 
  \usepackage{soul}

   \def\R{\mathbb{R}}
   \def\N{\mathbb{N}}
   \def\Z{\mathbb{Z}}
   \def\1{{\rm I\mskip -10.5mu 1}} 
   
   \def\e{{\varepsilon}}
   \def\D{{\nabla}}

   \def\cC{{\cal C}}

   \def\cF{{\cal F}}

   \def\cL{{\cal L}}

   \def\meas{\mathop{\rm meas}\nolimits}

   \def\Lip{\mathop{\rm Lip}\nolimits}
 
   \def\no{\noindent}
   \def\proof{\mbox {{\underline {\sf Proof}} \hspace{2mm}}}
   \def\qed{{\hfill {\em q.e.d.}\\\vspace{1mm}}}

   \newcommand{\beq}{\begin{equation}}
   \newcommand{\eeq}{\end{equation}}

\newtheorem{df}{Definition}[section]
\newtheorem{prop}[df]{Proposition}
\newtheorem{lemma}[df]{Lemma}
\newtheorem{teo}[df]{Theorem}
\newtheorem{rem}[df]{Remark}

\newtheorem{cor}[df]{Corollary}

 \newcommand{\sezione}[1]{\section{#1}\setcounter{equation}{0}}

\begin{document}


   \title{Infinitely many solutions for elliptic equations with
     non-symmetric nonlinearities
}

 \date{}

  \maketitle
\vspace{-2cm}



\begin{center}

{ {\bf Riccardo MOLLE$^a$,\quad Donato PASSASEO$^b$

\vspace{2mm}

\today
}}

\vspace{5mm}

{\em
${\phantom{1}}^a$Dipartimento di Matematica,
Universit\`a di Roma ``Tor Vergata'',\linebreak
Via della Ricerca Scientifica n. 1,
00133 Roma, Italy.}

\vspace{2mm}

{\em
${\phantom{1}}^b$Dipartimento di Matematica ``E. De Giorgi'',
  Universit\`a di Lecce,\linebreak 
P.O. Box 193, 73100 Lecce, Italy.
}
\end{center}

\vspace{5mm}


{\small {\sc \noindent \ \ Abstract.} - We deal with the existence of
  infinitely many solutions for a class of elliptic problems  
with non-symmetric nonlinearities.
Our result, which is motivated by a well known conjecture
formulated by A. Bahri and P.L. Lions, suggests a new approach to
tackle these problems. 
The proof is based on a method which does not require to use
techniques of deformation from the symmetry and may  
be applied to more general non-symmetric problems.
\footnote{ {\em E-mail address:} molle@mat.uniroma2.it (R. Molle).}
 
\vspace{3mm}



 \vspace{1mm}


}


\sezione{Introduction}


Let us consider the problem
\beq
\label{P}
-\Delta u=|u|^{p-1}u+w\quad\mbox{ in }\Omega,\qquad u=0\quad\mbox{ on
}\partial\Omega
\eeq
where $\Omega$ is a smooth bounded domain of $\R^n$, with $n\ge 1$,
$w\in L^2(\Omega)$, $p>1$ and $p<\frac{n+2}{n-2}$ when $n\ge 3$.

\no
If $w\not\equiv 0$ in $\Omega$, the corresponding energy functional
$E:H^1_0(\Omega)\to\R$, defined by 
\beq
E(u)=\frac12\int_\Omega|\D u|^2dx-\frac{1}{p+1}\int_\Omega
|u|^{p+1}dx-\int_\Omega w\, u\, dx
\eeq
is not even, so the equivariant Lusternik-Schnirelmann theory for
$\Z_2$-symmetric sets cannot 
be applied to find infinitely many solutions as in the case $w\equiv
0$ (see for instance \cite{A73,AR73,C69,C70,H70,H71,R74,R86} and also \cite{BC96,CP91} for a more
general framework). 

\no In the case $w\not\equiv 0$ in $\Omega$, a natural question (which
goes back to the beginning of the eighties) is wether the infinite
number of solutions still persists under perturbation. 

\no A detailed analysis was originally carried on in
\cite{A74,AR73,BB80,BB81,BB84,BL88,EG98,K64,MP75,R82,S80,T89} by Ambrosetti, Bahri, Berestycki,  Ekeland, Ghoussoub, Krasnoselskii,
Lions, Marino, Prodi, Rabinowitz, Struwe and Tanaka by introducing new perturbation methods. 
In particular, this question was raised to the attention  by
Rabinowitz also in his monograph on minimax methods (see \cite[Remark
10.58]{R86}).

In [2] Bahri proved that, if $n\ge 3$ and $1<p<\frac{n}{n-2}$, 
then there exists an open dense set of $w$ in $L^2(\Omega)$ such that
problem (\ref{P}) admits infinitely many solutions.
In [6] Bahri and Lions 
proved that, if $n\ge 3$ and $1<p<\frac{n}{n-2}$, then problem (\ref{P}) 
admits infinitely many solutions for every $w\in L^2(\Omega)$.
 
These results suggest the following conjecture, proposed by Bahri and
Lions in \cite{BL88}:  the multiplicity result obtained in \cite{BL88} holds also under the more general assumption $1<p<\frac{n+2}{n-2}$.
 
\no More recently, a new approach to tackle the break of symmetry in
elliptic problems has been developed by  Bolle,  Chambers, Ghoussoub
and Tehrani (see \cite{B99,BGT00,CG02}, which include also
applications to more general nonlinear problems).
However that approach did not allow to solve the Bahri-Lions
conjecture.

In the present paper we describe a new possible method to approach this problem.
By minimizing the energy functional $E$ in suitable subsets of $H^1_0(\Omega)$, we obtain infinitely many functions that present an arbitrarily large number of nodal regions having a prescribed structure (a check structure).
Their energy tends to infinity as the number of nodal regions tends to infinity.
Moreover, these functions satisfy equation \eqref{P} in each nodal region when the number of nodal regions is large enough (see Proposition \ref{Pn2.4}) and they are solutions of problem \eqref{P} when, in addition, they satisfy the assumptions of Proposition \ref{Pn2.5}. 

The idea is to trying to piece together solutions of Dirichlet problems in subdomains of $\Omega$ 
chosen in a suitable way.
This idea has been first used by Struwe in earlier papers (see \cite{S80,S81,S82} and references therein).
In the present paper we consider as nodal regions subdomains of $\Omega$ that are suitable deformations of cubes.
When the sizes of these cubes are all small enough, the nodal functions with check structure that we obtain seem to present suitable stability properties so that they persist when the problem \eqref{P} is perturbed by the term $w$.
The deformations of the nodal regions we use to construct solutions of problem \eqref{P}. are obtained in the present paper by considering a class of Lipshitz maps.
It is interesting to observe that such a class also appeared in some recent works of Rabinowitz and Byeon (see \cite{BR13,BR16} and the references therein) concerning a rather different problem: construct solutions having certain prescribed patterns for an Allen-Cahn model equation.
Also in that papers, as in the present one, Lipschitz condition is combined with the structure of $\Z^n$ and the covering of $\R^n$ by cubes with vertices in $\Z^n$.

In order to verify the assumptions of Proposition \ref{Pn2.5}, we need a technical condition (condition \eqref{2.20}).
In Lemma \ref{ML} we show that this condition is satisfied, for example, in the case $n=1$ (the proof may be also adapted to deal with radial solutions in domains $\Omega$ having radial symmetry).

Indeed, in dimension $n=1$, a more general result was obtained by Ehrmann in \cite{E57} (see also \cite{FL75,E57b} for related results).
Here it is proved that the ordinary differential equation
\beq
-u''(x)=f(u(x))+w(x)\qquad\text{ for }x\in (0,1),\ u(0)=u(1)=0
\eeq
has infinitely many distinct solutions when $f$ is a function with superlinear growth satisfying quite general assumptions.
However, the method here used relies on a shooting argument, typical of ordinary differential equations, combined with counting the oscillations of the solutions in the interval $(0,1)$.
Therefore, this method, which gives the existence of solutions having a sufficiently large number of zeroes in dimension $n=1$,  cannot be extended to higher dimensions.

On the contrary, in the present paper we use a method which is more similar to the one introduced by Nehari in \cite{N61}, that can be in a natural way extended to the case $n>1$.
In fact, for example, Nehari's work was followed up by Coffman who studied an analogous problem for partial differential equations (see \cite{C69,C70}).
Independently, this problem was also studied by Hempel (see \cite{H70,H71}).

More recently, the method introduced by Nehari has been also used by Conti, Terracini and Verzini to study optimal partition problems in $n$-dimensional domains and related problems: in particular, existence of minimal partitions and extremality conditions, behaviour of competing species systems with large interactions, existence of changing sign solutions for superlinear  elliptic equations, etc. (see \cite{CTV02,CTV03,CTV05,TV00}).

Notice that Nehari's work deal with an odd differential operator, so the corresponding energy functional is even.
Moreover, Nehari proves that for every positive integer $k$ there exists a solution having exactly $k$ zeroes.
On the contrary, in the present paper (as Ehrmann in \cite{E57}) we find only solutions with a large number of zeroes; moreover, we prove that, for all $w$ in $L^2(\Omega)$, the zeroes tend to be uniformly distributed in all of the domain $\Omega$ as their number tends to infinity (see Lemmas \ref{ML} and \ref{Ln3.3})
The reason is that, when $w\not\equiv 0$, the Nehari type argument we use in the proof works only when the sizes of all the nodal regions are small enough, so their number is sufficiently large.

In order to show that our existence result is sharp, we prove also that the term $w$ in problem \eqref{P} can be chosen in such a way that the problem does not have solutions with a small number of nodal regions.
More precisely, in the case $n=1$ we show that for all positive integer $h$ there exists $w_h$ in $L^2(\Omega)$ such that every solution of problem \eqref{P} with $w=w_h$ has at least $h$ zeroes (see Corollary \ref{Cr3}).
Indeed, we show that for all $n\ge 1$ and for every eigenfunction $e_k$ of the Laplace operator $-\Delta$ in $H^1_0(\Omega)$ there exists $\bar w_k$ in $L^2(\Omega)$ such that every solution $u$ of problem \eqref{P} with $w=\bar w_k$ must have the sign related to the sign of $e_k$ in the sense that every nodal region of $e_k$ has a subset where $u$ and $e_k$ have the same sign (see Proposition \ref{Pr2}). 

In the case $n>1$, condition (\ref{2.20}) seems to be more difficult to be verified because the class of all Lipschitz deformations of the nodal regions might result too large, as we explain in Remark \ref{Rn3.1}.
Therefore, in this case a useful idea might be to restrict the class of the admissible deformations so that we can verify a condition analogous to \eqref{2.20} and apply our method to construct nodal solutions having check structure.
For example, as we describe in Remark \ref{Rn3.1}, we can fix a suitable Lipschitz map $T_0:\overline\Omega\to\overline\Omega$ and consider nodal regions deformed by Lipschitz maps suitably close to $T_0$.
It is clear that, in order to apply our method, we need now to verify a condition analogous to \eqref{2.20} (that is condition \eqref{R*}) which holds or fails depending on the choice of $T_0$ and of the neighborhood of deformations close to $T_0$.
In a similar way, for example, we prove that if $\Omega$  is a cube of $\R^n$ with $n>1$, $p>1$, $p<\frac{n+2}{n-2}$ if $n>2$, for all $w$ in $L^2(\Omega)$ there exist infinitely many solutions $u_k(x)$ of problem \eqref{P} such that the nodal regions of the function $u_k\left(\frac xk\right)$, after translations, tend to the cube as $k\to\infty$ (the proof will be reported in a paper in preparation).
We believe that this result may be extended to every interval or pluri-interval of $\R^n$ with $n>1$ and then to every bounded domain $\Omega$ by a suitable choice of the deformation $T_0$, related to the geometrical properties of the domain $\Omega$. 

Let us point out that our method does not require techniques of
deformation from the symmetry and may be applied to more general
problems: for example, when the nonlinear term $|u|^{p-1}u$ is
replaced by $c_+(u^+)^p-c_-(u^-)^p$ with $c_+$ and $c_-$ two positive
constants (see Lemma \ref{Ln3.3}), in case of different, nonhomogeneous boundary conditions
and even in case of nonlinear elliptic equations involving critical
Sobolev exponents.

\vspace{2mm}

{\small {\bf Acknowledgement}. The authors are very much grateful to Professor P.H. Rabinowitz for several helpful comments, suggestions and informations on this work and on related subjects.
}


\sezione{Existence of infinitely many nodal solutions} 
\label{S2}


In order to find infinitely many solutions with an arbitrarily large
number of nodal regions, we proceed as follows.

\no Let us set 
\beq
\begin{array}{c}
\vspace{2mm}
C_0=\{x\in\R^n\ :\ 0<x_i<1\ \mbox{ for }i=1,\ldots,n\},\\
C_z=z+C_0,\quad
\sigma(z)=(-1)^{\sum_{i=1}^n z_i}\qquad\forall
z\in\Z^n,
\end{array}
\eeq
\beq
Z_k=\left\{z\in\Z^n\ :\ {1\over k}\, C_z\subset\Omega\right\},\qquad
P_k=\bigcup_{z\in Z_k}{1\over k}\, \overline{C}_z,\qquad \forall k\in\N.
\eeq
Notice that there exists $k_\Omega$ in $\N$ such that
$Z_k\neq\emptyset$ $\forall k\ge k_\Omega$.

\no For all subsets $P,Q$ of $\R^n$ and for all $L\ge 1$, let us
denote by $\cC_L(P,Q)$ the set of all the functions $T:P\to Q$ such
that 
\beq
{1\over L}\,|x-y|\le |T(x)-T(y)|\le L\, |x-y|\qquad\forall x,y\in P.
\eeq

\no For all $k\ge k_\Omega$, $z\in Z_k$, $L\ge 1$, $T\in
\cC_L(P_k,\overline\Omega)$ let us set
\beq
\label{1.1}
E^T_{k,z}=\inf\left\{E(u)\ :\ u\in H^1_0\left(T\left({1\over k}\,
      C_z\right)\right),\ \int_{T\left({1\over k}\,
      C_z\right)}|u|^{p+1}dx=1\right\}.
\eeq
Since $p<{n+2\over n-2}$ when $n\ge 3$, one can easily verify that the
infimum in (\ref{1.1}) is achieved.
Moreover, for all $L\ge 1$ and $k\ge k_\Omega$, also the infimum 
\beq
\inf\{E^T_{k,z}\ :\ z\in Z_k,\ T\in C_L(P_k,\overline\Omega)\}
\eeq
is achieved (as one can prove by standard arguments using
Ascoli-Arzel\`a Theorem) and the following lemma holds.
\begin{lemma}
\label{Ln2.2}
For all $L\ge 1$, we have
\beq
\lim_{k\to\infty}\min\{E ^T_{k,z}\ :\ z\in Z_k,\
T\in\cC_L(P_k,\overline\Omega)\}=\infty.
\eeq
\end{lemma}
\no \proof 
For all $L\ge 1$ and $k\ge k_\Omega$, let us choose $z_k\in Z_k$,
$T_k\in \cC_L(P_k,\overline\Omega)$ and $\bar u_k\in
H^1_0\left(T_k\left({1\over k}C_{z_k}\right)\right)$
such that 
\beq
\label{1.2}
\int_{T_k\left({1\over k}C_{z_k}\right) }|\bar u_k|^{p+1}dx=1
\ \mbox{  and }\
E(\bar u_k)=E^{T_k}_{k,z_k}=\min\{E^T_{k,z}\ :\ z\in Z_k,\ T\in
\cC_L(P_k,\overline\Omega)\}.
\eeq
We say that 
\beq
\label{Finfty}
\lim_{k\to\infty}\int_{T_k({1\over
  k}\,C_{z_k})}|\D \bar u_k|^{2}dx=\infty.
\eeq
In fact, arguing by contradiction, assume that
\beq
\liminf_{k\to\infty}\int_{T_k({1\over
  k}\,C_{z_k})}|\D\bar  u_k|^{2}dx<\infty.
\eeq
It follows that (up to a subsequence)  $(\bar u_k)_k$ is bounded in
$H^1_0(\Omega)$ and there exists a function $\bar u\in H^1_0(\Omega)$ such
that $\bar u_k\to \bar u$, as $k\to\infty$, weakly in $H^1_0(\Omega)$, in
$L^{p+1}(\Omega)$, and almost everywhere in $\Omega$ (here $\bar u_k$ is
extended by the value 0 in $\Omega\setminus T_k({1\over
  k}\,C_{z_k})$).
Since $\meas\left(T_k({1\over k}\,C_{z_k})\right)\to 0$ as $k\to\infty$, from
  the almost everywhere convergence we obtain $\bar u\equiv 0$ in
  $\Omega$, which is a contradiction because $\bar u_k\to\bar u$
  in $L^{p+1}(\Omega)$ and (\ref{1.2}) holds for all $k\ge k_\Omega$.
Thus (\ref{Finfty}) is proved.

\no Notice that
\beq
E^{T_k}_{k,z_k}={1\over 2}\int_{T_k\left({1\over k}\,
    C_{z_k}\right)}|\D\bar u_k|^2dx-{1\over p+1}-
\int_{T_k\left({1\over k}\,
    C_{z_k}\right)} \bar u_kw\, dx\qquad\forall k\ge k_\Omega
\eeq
where
\beq
\begin{array}{rcl}
\left| 
\int_{T_k({1\over k}\,C_{z_k})} \bar u_k\, w\, dx
\right| 
& \le &
\left(\int_{T_k({1\over k}\,C_{z_k})}\bar  u_k^2 dx\right)^{1\over 2}
\left(\int_{T_k({1\over k}\,C_{z_k})}  w^2 dx\right)^{1\over 2}
\\ [2ex]
& \le &
\left[\meas\left(T_k({1\over k}\,C_{z_k})\right)\right]^{{1\over
    2}-{1\over p+1}}\left(\int_{\Omega}  w^2 dx\right)^{1\over 2}.
\end{array}
\eeq
As a consequence, for all  $k\ge k_\Omega$ we obtain
\beq
E^{T_k}_{k,z_k}\ge {1\over 2}
\int_{T_k({1\over k}\,C_{z_k})}|\D \bar u_k|^2 dx-{1\over p+1}-
\left(\int_{\Omega}  w^2 dx\right)^{1\over 2}\cdot 
\left[\meas\left(T_k\left({1\over k}\,C_{z_k}\right)\right)\right]^{{1\over
    2}-{1\over p+1}},
\eeq
and, as $k\to\infty$,
\beq
\lim_{k\to\infty}E^{T_k}_{k,z_k}=\infty 
\eeq
which completes the proof.

\qed

\begin{cor}
\label{C2.2}
For all $L\ge 1$ there exists $k(L)\ge k_\Omega$ such that for all
$k\ge k(L)$, $z\in Z_k$ and $T\in \cC_L(P_k,\overline\Omega)$ the
minimum 
\beq
\label{2.3}
\min\left\{E(u)\ :\ u\in H^1_0\left( T\left({1\over k}\,
      C_z\right)\right),\ \int_{T\left({1\over k}\,
      C_z\right) }|u|^{p+1}dx<1\right\}
\eeq
is achieved by a unique minimizing function $\tilde u^T_{k,z}$.
Moreover, we have
\beq
\label{2.4}
\lim_{k\to\infty}\sup\left\{ \int_{T\left({1\over k}\,
      C_z\right) }|\D\tilde u^T_{k,z}|^2dx\ :\ z\in Z_k,\
  T\in\cC_L(P_k,\overline \Omega)\right\}=0.
\eeq
\end{cor}
\no \proof As a consequence of Lemma \ref{Ln2.2}, for all $L\ge 1$
there exists $k(L)\ge k_\Omega$ such that 
$$
0<\min\left\{E(u)\ :\ u\in H^1_0\left(T\left({1\over k}\,
      C_z\right)\right),\ \int_{T\left({1\over k}\,
      C_z\right) }|u|^{p+1}dx=1\right\}
$$
\beq
\hspace{3cm} \forall k\ge k(L),\
\forall z\in Z_k,\ \forall T\in \cC_L(P_k,\overline\Omega).
\eeq
On the other hand,
$$
\inf\left\{E(u)\ :\ u\in H^1_0\left(T\left({1\over k}\,
      C_z\right)\right),\ \int_{T\left({1\over k}\,
      C_z\right)}|u|^{p+1}dx<1\right\}\le 0
$$
\beq
\label{inf}\hspace{3cm} \forall k\ge k_\Omega,\
\forall z\in Z_k,\ \forall T\in \cC_L(P_k,\overline\Omega)
\eeq
because $E(u)=0$ for $u\equiv 0$ in $T\left({1\over k}\,
      C_z\right)$.

\no Now, let us consider a minimizing sequence for the infimum in
(\ref{inf}).
Since it is bounded in $L^{p+1}\left({1\over k}\, C_z\right)$, we
infer from (\ref{inf}) that it is bounded also in $H^1_0\left({1\over
    k}\, C_z\right)$. 
Therefore, since $p<{n+2\over n-2}$ when $n\ge 3$, one can prove by
standard arguments that (up to a subsequence) it converges to a
function $\tilde u^T_{k,z}\in H^1_0\left(T\left({1\over k}\, C_z\right)\right)$ such
that $\int_{T\left({1\over k}\, C_z\right)}|\tilde
u^T_{k,z}|^{p+1}dx$ $<1$ and
\beq
E(\tilde u^T_{k,z})=\min\left\{E(u)\ :\ u\in H^1_0\left(T\left({1\over k}\,
  C_z\right)\right),\ \int_{T\left({1\over k}\,
      C_z\right)}|u|^{p+1}dx<1\right\}.
\eeq
In order to prove (\ref{2.4}) we argue by contradiction and assume
that 
\beq
\limsup_{k\to\infty}\sup\left\{ 
\int_{T\left({1\over k}\, C_z\right)} |\D\tilde u^T_{k,z}|^2dx\ :\ z\in
Z_k,\ T\in \cC_L(P_k,\overline\Omega)\right\}>0.
\eeq
Then, for all $k\ge k(L)$ there exist $z_k\in Z_k$ and
$T_k\in\cC_L(P_k,\overline\Omega)$ such that (up to a subsequence)
\beq
\label{du>0}
\lim_{k\to\infty} 
\int_{T_k\left({1\over k}\, C_{z_k}\right)} |\D\tilde u^{T_k}_{k,z_k}|^2dx>0.
\eeq
Since $E(\tilde u^{T_k}_{k,z_k})\le 0$ and the sequence $\tilde
u^{T_k}_{k,z_k}$ (extended by the value zero outside ${1\over k}\,
C_{z_k}$) is bounded in $L^{p+1}(\Omega)$, we infer that it is bounded
also in $H^1_0(\Omega)$.
We say that, as a consequence, $\tilde u^{T_k}_{k,z_k}\to 0$ as
$k\to\infty$ in $L^{p+1}(\Omega)$.
In fact, since the sequence $\left(\tilde
  u^{T_k}_{k,z_k}\right)_{k\in\N}$ is bounded in $H^1_0(\Omega)$, it
converges weakly in $H^1_0(\Omega)$, in $L^{p+1}(\Omega)$ and a.e. in
$\Omega$ to a function $\tilde u\in H^1_0(\Omega)$.
Since $\lim\limits_{k\to\infty}\meas \left({1\over k}\,
  C_{z_k}\right)=0$, we can say that $\tilde u\equiv 0$ in $\Omega$.
Thus, $\tilde u^{T_k}_{k,z_k}\to 0$ as $k\to\infty$ in
$L^{p+1}(\Omega)$.
Therefore, taking into account that
\beq
E\left(\tilde u^{T_k}_{k,z_k}\right)={1\over 2}\int_\Omega |\D 
\tilde u^{T_k}_{k,z_k}|^2dx-{1\over p+1} \int_\Omega |\tilde
u^{T_k}_{k,z_k}|^{p+1} dx -  \int_\Omega w \, \tilde
u^{T_k}_{k,z_k}\, dx\le 0
\eeq
it follows that $\tilde u^{T_k}_{k,z_k}\to 0$ also in $H^1_0(\Omega)$
in contradiction with (\ref{du>0}).

\no Thus, we can conclude that (\ref{2.4}) holds.
Finally, notice that $\tilde u^{T}_{k,z}$ is the unique minimizing
function for (\ref{2.3}) because the functional $E$ is strictly convex
in a suitable neighborhood of zero.
So the proof is complete.

\qed

Taking into account Corollary \ref{C2.2}, for all $k\ge k(L)$, $z\in
Z_k$ and $T\in\cC_L(P_k,\overline \Omega)$ we can consider a
minimizing function $\tilde u^T_{k,z}$ for the minimum (\ref{2.3}).
Moreover, since $p>1$, for all $u\in H^1_0
\left({1\over k}\, C_{z}\right)$ there exists the maximum
\beq
M(u)=\max\left\{E\left(\tilde u^T_{k,z}+t(u-\tilde u^T_{k,z})\right)\ :\
t\ge 0\right\}
\eeq
and $M(u)\ge E^T_{k,z}$ when $u\not\equiv \tilde u^T_{k,z}$ in
$T\left({1\over k}\, C_{z}\right)$.

\begin{lemma}
\label{Ln2.3}
For all $k\ge k(L)$, $z\in Z_k$ and $T\in\cC_L(P_k,\overline\Omega)$,
there exists a function $u^T_{k,z}$ in $ H^1_0
\left(T\left({1\over k}\, C_{z}\right)\right)$ such that
$u^T_{k,z}\not\equiv \tilde u^T_{k,z}$, $\sigma(z)[u^T_{k,z}-\tilde
u^T_{k,z}]\ge 0$ in $T\left({1\over k}\, C_{z}\right)$ and 
$$
\hspace{-2cm}
E(u^T_{k,z})=M(u^T_{k,z})=\min
\left\{M(u)\ :\ u\in H^1_0\left( T\left({1\over k}\,
      C_{z}\right)\right),\ u\not\equiv \tilde u^T_{k,z},\right. 
$$
\beq
\label{min}
\hspace{5cm}\left.{\phantom{****{1\over k}}} \sigma(z)[u-\tilde
u^T_{k,z}]\ge 0\ \mbox{ in }\ T\left({1\over k}\,
  C_{z}\right)\right\}.
\eeq
Moreover, we have $E(u^T_{k,z})\ge E^T_{k,z}$.
\end{lemma}

\no \proof Let us consider a minimizing sequence $(u_i)_{i\in\N}$ for
the minimum (\ref{min}).
Whitout any loss of generality, we can assume that
\beq
\label{p+1}
\int_{T\left({1\over k}\, C_z\right) }|u_i-\tilde u^T_{k,z}|^{p+1}dx
=1\qquad\forall i\in\N.
\eeq
It follows that this sequence is bounded in $H^1_0\left(T\left({1\over
      k}\, C_z\right)\right)$.
Therefore, since $p<{n+2\over n-2}$ when $n\ge 3$, up to a subsequence
it converges weakly in $H^1_0$, in $L^{p+1}$ and a.e. to a function
$\hat u\in H^1_0\left(T\left({1\over
      k}\, C_z\right)\right)$.

\no Notice  that the $L^{p+1}$ convergence and (\ref{p+1}) imply
\beq
\int_{T\left({1\over k}\, C_z\right)} |\hat u-\tilde u^T_{k,z}|^{p+1}dx=1,
\eeq
so $\hat u\not\equiv \tilde u^T_{k,z}$.
Moreover, the a.e. convergence implies $\sigma(z)[\hat u-\tilde
u^T_{k,z}]\ge 0$ in $T\left({1\over k}\, C_z\right)$. 
We say that, indeed, the convergence is strong in $H^1_0\left(
  T\left({1\over k}\, C_z\right)\right)$. 
In fact, arguing by contradiction, assume that (up to a subsequence)
\beq
\label{<}
\int_{T\left({1\over k}\, C_z\right)} |\D\hat u |^2dx<\lim_{i\to\infty}
\int_{T\left({1\over k}\, C_z\right)} |\D u_i |^2dx.
\eeq
As a consequence, we obtain $M(\hat u)<\lim\limits_{i\to\infty}M(u_i)$
which is a contradiction because $\hat u\not\equiv \tilde u^T_{k,z}$
and $(u_i)_{i\in\N}$ is a minimizing sequence for (\ref{min}) so that
$M(\hat u)\ge \lim\limits_{i\to\infty}M(u_i)$. 

\no Therefore, we can conclude that $u_i\to \hat u$ in $H^1_0\left(
  T\left({1\over k}\, C_z\right)\right)$ and
$\lim\limits_{i\to\infty}M(u_i)=M(\hat u)$. 
Since $p>1$, there exists $\hat t>0$ such that 
\beq
E(\tilde u^T_{k,z}+\hat t(\hat u-\tilde u^T_{k,z}))=M(\hat u).
\eeq
Thus, all the assertions of Lemma \ref{Ln2.3} hold with $ u^T_{k,z} =
\tilde u^T_{k,z}+\hat t(\hat u-\tilde u^T_{k,z})$. 

\qed

\begin{prop}
\label{Pn2.4}
There exists $k_1(L)\ge k(L)$ such that for all $k\ge k_1(L)$, $z\in
Z_k$ and $T\in\cC_L(P_k,\overline\Omega)$ the function $u^T_{k,z}$ is a
solution of the Dirichlet problem  
\beq
\label{eD}
-\Delta u=|u|^{p-1}u+w\quad{\rm in }\ T\left({1\over k}\,
  C_z\right),\qquad u=0\quad{\rm on }\ \partial T\left({1\over k}\,
  C_z\right). 
\eeq
\end{prop}

\no\proof It is clear that the function $\tilde u^T_{k,z}$ (local
minimum of the functional $E$) is a solution of the Dirichlet problem
(\ref{eD}). 
In order to prove that, for $k$ large enough,  also $ u^T_{k,z}$ is
solution of the same problem, let us consider the function
$G:T\left({1\over k}\, C_z\right)\times\R\to \R$ such that
$G(x,\cdot)\in \cC^2(\R)$  $\forall x\in T\left({1\over k}\, C_z\right)$
and 
\beq
\begin{array}{cl}
\vspace{3mm}G(x,t) = 
{|t|^{p+1}\over p+1}+w(x)\, t &{\rm if }\ \sigma(z)[t-\tilde u^T_{k,z}(x)]\ge 0,
\\
{\partial^2  G\over \partial t^2}\, (x,t)={\partial^2  G\over \partial t^2}\,
(x,\tilde u^T_{k,z}(x)) &{\rm if }\ \sigma(z)[t-\tilde u^T_{k,z}(x)]\le 0.
\end{array}
\eeq
Moreover, let us set $g(x,t)= \frac{dG}{dt}(x,t)$.

\noindent Then, consider the functional $E_{k,z,T}:H^1_0\left(T
  \left({1\over k}\, C_z\right)\right)\to\R$ defined by 
\beq
E_{k,z,T}(u)={1\over 2}\int_{T\left({1\over k}\, C_z\right)}|\D u|^2dx-
\int_{T\left({1\over k}\, C_z\right)}G(x,u)\, dx.
\eeq

\no Let us assume, for example, $\sigma(z)=1$ (in a similar way one can
argue when $\sigma(z)=-1$).

\no One can verify by direct computation that for all $u\not\equiv\tilde
u^T_{k,z}$ there exists a unique $t_u>0$ such that 
\beq
E'_{k,z,T}(\tilde u^T_{k,z}+t_u(u-\tilde u^T_{k,z}))[u-\tilde
u^T_{k,z}]=0
\eeq
if and only if $(u-\tilde u^T_{k,z})\vee 0\not\equiv 0$.
In this case $ E'_{k,z,T}(\tilde u^T_{k,z}+t(u-\tilde u^T_{k,z})[u-\tilde
u^T_{k,z}]$ is positive for $t\in ]0,t_u[$ and negative for $t>t_u$,
so we have 
\beq
E_{k,z,T}(\tilde u^T_{k,z}+t_u(u-\tilde u^T_{k,z}))
=
\max\{
E_{k,z,T}(\tilde u^T_{k,z}+t(u-\tilde u^T_{k,z}))\ :\ t >0\}.
\eeq
Moreover, we have
\beq
E''_{k,z,T}(\tilde u^T_{k,z}+t_u(u-\tilde u^T_{k,z}))[u-\tilde
u^T_{k,z},u-\tilde u^T_{k,z}]<0.
\eeq
Taking into account that $\tilde u^T_{k,z}$ is solution of problem
(\ref{eD}), we obtain by direct computation 
\begin{eqnarray}
\nonumber
E_{k,z,T}(\tilde u^T_{k,z}+t(u-\tilde u^T_{k,z}))
& = &
E_{k,z,T}(\tilde u^T_{k,z}+t[(u-\tilde u^T_{k,z})\vee 0])
\\
& &
+{t^2\over 2} \int_{T\left({1\over k}\, C_z\right)}
|\D [(u-\tilde u^T_{k,z})\wedge 0] |^2dx \\
\nonumber
& &
- {t^2\over 2}\, p  \int_{T\left({1\over k}\, C_z\right)}
|\tilde u^T_{k,z}|^{p-1}[(u-\tilde u^T_{k,z})\wedge 0]^2dx
\qquad\forall t>0.
\end{eqnarray}
Notice that
$$
\hspace{-6cm}
\int_{T\left({1\over k}\, C_z\right)}
|\tilde u^T_{k,z}|^{p-1}[(u-\tilde u^T_{k,z})\wedge 0]^2dx
$$
\beq
\hspace{3cm}
\le
\left(\int_{T\left({1\over k}\, C_z\right)}
|\tilde u^T_{k,z}|^{p+1}dx\right)^{p-1\over p+1}
\left(\int_{T\left({1\over k}\, C_z\right)}
|(u-\tilde u^T_{k,z})\wedge 0|^{p+1}dx\right)^{2\over p+1}
\eeq
and, by (\ref{2.4}), 
\beq
\lim_{k\to\infty}\sup
\left\{ \int_{T\left({1\over k}\, C_z\right)}|\tilde u^T_{k,z}|^{p+1}dx
\ :\
z\in Z_k,\ T\in\cC_L(P_k,\overline\Omega)\right\}=0.
\eeq
Moreover, we have
$$
\int_{T\left({1\over k}\, C_z\right)}|\D [(u-\tilde u^T_{k,z})\wedge 0]
|^2dx
\ge
\lambda_k \left(\int_{T\left({1\over k}\, C_z\right)}
|(u-\tilde u^T_{k,z})\wedge 0|^{p+1}dx\right)^{2\over p+1}
$$
\beq
\hspace{3cm}
\forall k\in \N,\ \forall z\in Z_k,\ \forall
T\in\cC_L(P_k,\overline\Omega)
\eeq
where, for all $k\in\N$,
$$
\hspace{-5cm} \lambda_k=\inf\left\{\int_{T\left({1\over k}\,
      C_z\right)} |\D\psi|^2dx\ :\ z\in Z_k,\ T\in
  \cC_L(P_k,\overline\Omega),\ \right. 
$$
\beq
\hspace{5cm} \left. \psi\in H^1_0\left(T\left({1\over k}\,
      C_z\right)\right),\ \int_{T\left({1\over k}\,
      C_z\right)}|\psi|^{p+1}dx=1\right\}. 
\eeq
Notice that $\lim\limits_{k\to\infty}\lambda_k=\infty$ otherwise for
all $i\in\N$ there exist $k_i\in\N$ $z_i\in Z_{k_i}$,
$T_i\in\cC_L(P_{k_i},\overline\Omega)$, $\psi_i\in
H^1_0\left(T_i\left({1\over k_i}\, C_{z_i}\right)\right)$ (extended by
the value zero outside $T_i\left({1\over k_i}\, C_{z_i}\right)$) such
that $\lim\limits_{i\to\infty} k_i=\infty$,
$\int_{\Omega}|\psi_i|^{p+1}dx=1$ $\forall i\in\N$ and
$\lim\limits_{i\to\infty}\int_\Omega|D\psi_i|^2dx<\infty$. 

\no As a consequence, since $p<{n+2\over n-2}$ when $n\ge 3$, there
exists $\bar\psi$ in $H^1_0(\Omega)$ such that (up to a subsequence)
$\psi_i\to \bar\psi$ as $i\to\infty$ weakly in $H^1_0(\Omega)$, in
$L^{p+1}(\Omega)$ and a.e. in $\Omega$. 

\no Moreover, since $\lim\limits_{i\to\infty}\meas T_i\left({1\over
    k_i}\, C_{z_i}\right)=0$, the a.e. convergence implies that
$\bar\psi\equiv 0$ in $\Omega$, which is in contradiction with the
convergence in $L^{p+1}(\Omega)$ because
$\int_{\Omega}|\psi|^{p+1}dx=1$ $\forall i\in\N$. 

\no Thus, we can conclude that $\lim\limits_{k\to\infty}\lambda_k=\infty$.
 It follows that, for $k$ large enough,
\beq
E_{k,z,T}(\tilde u^T_{k,z}+t[(u-\tilde u^T_{k,z})\vee 0])
\le
E_{k,z,T}(\tilde u^T_{k,z}+t(u-\tilde u^T_{k,z}))\qquad \forall t>0,\
\forall u\in H^1_0\left({1\over k}\, C_z\right).
\eeq
As a consequence, if we denote by $\Gamma$ the set defined by 
\beq
\Gamma=\left\{u\in  H^1_0\left( T \left({1\over k}\, C_z\right)\right)
\ :\ 
u\not\equiv \tilde u^T_{k,z},\ E_{k,z,T}'(u)[u-\tilde u^T_{k,z}]=0\right\},
\eeq
we have $u^T_{k,z}\in\Gamma$ and 
\beq
E(u^T_{k,z})= E_{k,z,T}(u^T_{k,z})=\min_\Gamma E_{k,z,T}.
\eeq
Therefore, there exists a Lagrange multiplier $\mu\in\R$ such that 
\beq
E'_{k,z,T}(u^T_{k,z})[\varphi]=
\mu\left\{E''_{k,z,T}(u^T_{k,z})[u^T_{k,z}-\tilde u^T_{k,z},\varphi]
+E'_{k,z,T}(u^T_{k,z})[\varphi]\right\}\quad \forall\varphi\in H^1_0
\left({1\over k}\, C_z\right).
\eeq
In particular, if we choose $\varphi=u^T_{k,z}-\tilde u^T_{k,z}$, we
obtain $\mu=0$ because $E'_{k,z,T}(u^T_{k,z})[u^T_{k,z}-\tilde
u^T_{k,z}]=0$ while $E''_{k,z,T}(u^T_{k,z})[u^T_{k,z}-\tilde u^T_{k,z},u^T_{k,z}-\tilde
u^T_{k,z}]\neq 0$.
Thus, $u^T_{k,z}$ is a weak solution of the Dirichlet problem
\beq
-\Delta u^T_{k,z}=g(x,u^T_{k,z})\quad{\rm in }\ T\left({1\over k}\,
  C_z\right),\qquad u=0\quad{\rm on }\ \partial T\left({1\over k}\,
  C_z\right).
\eeq
On the other hand, since $u^T_{k,z}-\tilde u^T_{k,z}\ge 0$ in $T\left({1\over k}\,
  C_z\right)$, we have
\beq
g(x,u^T_{k,z}(x)) =|u^T_{k,z}(x)|^{p-1}u^T_{k,z}(x)+w(x)\qquad \forall
x\in T\left({1\over k}\, C_z\right),
\eeq
so $u^T_{k,z}$ is a solution of problem (\ref{eD}).

\qed

\no When the function $u^T_k=\sum\limits_{z\in Z_k}u^T_{k,z}$
satisfies a suitable stationarity property, then it is solution of
problem (\ref{P}) (here the function $u^T_{k,z}$ is extended by the
value zero outside $T\left({1\over k}\, C_z\right)$).
In fact, the following proposition holds.

\begin{prop}
\label{Pn2.5}
Assume that $k\ge k_1(L)$ and $T\in \cC_L(P_k,\overline\Omega)$.
Moreover, assume that the function  $u^T_k=\sum\limits_{z\in
  Z_k}u^T_{k,z}$ satisfies the following condition: $E'(u^T_k)[v\cdot D
u_k^T]=0$ for all vector field $v\in\cC^1(\overline\Omega,\R^n)$ such
that $v\cdot\nu=0$ on $\partial\Omega$ (here $\nu$ denotes the outward
normal vector on $\partial\Omega$).
Then, $u_k^T$ is a solution of the Dirichlet problem (\ref{P}).
\end{prop}

\no \proof We have to prove that $E'(u^T_k)[\varphi]=0$ $\forall
\varphi\in H^1_0(\Omega)$.
Taking into account Proposition \ref{Pn2.4}, since $u^T_{k,z}$
satisfies the Dirichlet problem (\ref{eD}) for all $z\in Z_k$, we have 
\begin{eqnarray}
\nonumber
E'(u^T_k)[\varphi] & = & \int_\Omega [\D u^T_k\cdot \D\varphi -
|u^T_k|^{p-1}u^T_k\varphi-w\varphi]dx\\
& = &
\sum\limits_{z\in Z_k}\int_{T\left({1\over k}\, C_z\right)}
[\D u^T_k\cdot \D\varphi -
|u^T_k|^{p-1}u^T_k\varphi-w\varphi]dx\\
& = &
\sum\limits_{z\in Z_k} \int_{\partial T\left({1\over k}\, C_z\right)}
\varphi \, (\D u^T_k\cdot\nu_{k,z})d\sigma,
\end{eqnarray}
where $\nu_{k,z}$ denotes the outward normal on $\partial
T\left({1\over k}\, C_z\right)$.
Thus, in order to obtain $E'(u^T_k)$ $[\varphi]=0$, we have to prove that
if $z_1,z_2\in Z_k$ and $|z_1-z_2|=1$ (that is ${T\left({1\over k}\,
    C_{z_1}\right)}$ and ${T\left({1\over k}\, C_{z_2}\right)}$ are
adjacent subdomains of $\Omega$) then 
\beq
\label{nabla}
\D u^T_{k,z_1}(x)=\D u^T_{k,z_2}(x) \qquad \forall x\in\partial 
T\left({1\over k}\, C_{z_1}\right)\cap T\left({1\over k}\,
  C_{z_2}\right).
\eeq
Taking into account that $u^T_{k,z}$ satisfies problem (\ref{eD}) for
all $z\in Z_k$, for all vector field $v\in\cC^1(\overline\Omega,\R^n)$ such
that $v\cdot\nu=0$ on $\partial\Omega$ we obtain

$$
\hspace{-12cm} E'(u^T_k)[v\cdot  \D u^T_k]=
$$
 
\vspace{-8mm}

\begin{eqnarray}
\nonumber
 &= & 
\int_{\Omega} [\D u^T_k\cdot \D(v\cdot \D u^T_k)-|u^T_k|^{p-1}u^T_k 
(v\cdot \D u^T_k)-w(v\cdot \D u^T_k)]dx
\\
& = &
\nonumber
\sum\limits_{z\in Z_k}
\int_{T\left({1\over k}\, C_{z}\right)} [\D u^T_k\cdot \D(v\cdot
\D u^T_k)-|u^T_k|^{p-1}u^T_k (v\cdot \D u^T_k)-w(v\cdot \D u^T_k)]dx 
\\
& = & 
\sum\limits_{z\in Z_k}
\int_{\partial T\left({1\over k}\, C_{z}\right)} 
( \D u^T_k\cdot \nu_{k,z})^2(v\cdot \nu_{k,z})d\sigma.
\end{eqnarray}
Since $E'(u^T_k)[v\cdot \D u^T_k]=0$ $\forall v\in
\cC^1(\overline\Omega,\R^n)$ such that $v\cdot\nu=0$ on
$\partial\Omega$, (\ref{nabla}) follows easily.
Thus, we can conclude that $u^T_k$ is a solution of problem (\ref{P}).

\qed

\no In order to obtain a function $u^T_k$ which is stationary in the
sense of Proposition \ref{Pn2.5}, we can, for example, minimize $E(u^T_k)$ with respect
to $T$ for $k$ large enough.

\no First notice that, since $\Omega$ is a smooth bounded domain,
there exist $k'_\Omega\ge k_\Omega$ and $L'_\Omega\ge 1$ such that,
for all $k\ge k'_\Omega$ and $L\ge L_\Omega'$, we have 
\beq
\label{Fdiv0}
\{T\in\cC_L(P_k,\overline\Omega)\ :\
T(P_k)=\overline\Omega\}\neq\emptyset.
\eeq
Moreover, using Ascoli-Arzel\`a Theorem, one can show the following
lemma. 
\begin{lemma}
\label{L2.2}
If (\ref{Fdiv0}) holds, there exists $\widetilde
T_k^L\in\cC_L(P_k,\overline\Omega)$ such that $\widetilde
T^L_k(P_k)=\overline\Omega$ and
\beq
\sum_{z\in Z_k} E(u^{\widetilde T^L}_{k,z})=\min\left\{\sum_{z\in Z_k} E(u^{
  T}_{k,z})\ :\ T\in\cC_L(P_k,\overline\Omega),\ T(P_k)=\overline\Omega\right\}.
\eeq
\end{lemma}

\no For all $L\ge 1$ and $T\in\cC_L(P_k,\overline\Omega)$, let us set 
\beq
\cL(T)=\inf\left\{\cL\ :\ \cL\ge 1,\ {1\over\cL}\,
  |x-y|\le|T(x)-T(y)|\le\cL\, |x-y|\quad \forall x,y\in P_k\right\}.
\eeq
Using again Ascoli-Arzel\`a Theorem, we infer that, for all $L\ge
L'_\Omega$ and $k\ge k'_\Omega$,  there exists
$T_k^L\in\cC_L(P_k,\overline\Omega)$ such that $T_k^L\left({1\over k}\,
  C_z\right)=\widetilde T_k^L\left({1\over k}\,
    C_z\right)$ $\forall z\in Z_k$ and
\beq
\cL(T_k^L)=\min\left\{\cL(T)\ :\ T\in\cC_L(P_k,\overline\Omega),\
  T\left({1\over k}\, C_z\right)=\widetilde T_k^L\left({1\over k}\,
    C_z\right)\ \forall z\in Z_k\right\}.
\eeq
Notice that $T_k^L$ depends only on the geometrical properties of the
subdomains $\widetilde T_k^L\left({1\over k}\, C_z\right)$ with $z\in
Z_k$.
A large $\cL(T_k^L)$ means that there are large differences in the sizes
and in the shape of these subdomains.

\no We can now state the following multiplicity result.

\begin{teo}
\label{MT}
Let $n\ge 1$, $p>1$ and $p<{n+2\over n-2}$ when $n\ge 3$.
Moreover, assume that there exists $\bar L\ge  L'_\Omega$  such that 
\beq
\label{2.20}
\limsup_{k\to\infty}\cL(T_k^{\bar L})< \bar L.
\eeq
Then, problem (\ref{P}) admits infinitely many solutions
(see also Remark \ref{Rn3.1} concerning condition (\ref{2.20})).
\end{teo}

\no Theorem \ref{MT} is a direct consequence of the following
proposition.
\begin{prop}
\label{MP}
If the assumptions of Theorem \ref{MT} are satisfied, for all $w\in
L^2(\Omega)$ there exists  $\bar k\ge k_\Omega$ such
that for all $k\ge\bar k$ there exists $T^{\bar L}_k\in\cC_{\bar
  L}(P_k,\overline \Omega)$ satisfying the following property:
$T^{\bar L}_k(P_k)=\overline\Omega$ and the function
$u_k=\sum\limits_{z\in Z_k}u^{T^{\bar L}_k}_{k,z}$ is a solution of
problem (\ref{P}). 
Moreover, the number of nodal regions of $u_k$ tends to infinity as
$k\to \infty$ and 
\beq
\label{ktoinfty}
\lim_{k\to\infty}\min\left\{ E\left(u^{T^{\bar L}_k}_{k,z}\right)\ :\
  z\in Z_k\right\}=\infty.
\eeq
\end{prop}

\no\proof As a consequence of condition (\ref{2.20}), there exist
$\bar k\in\N$ and $\bar \e>0$ such that   
\beq
\label{D}
D\circ T_k^{\bar L}\in\cC_{\bar L}(P_k,\overline\Omega)\quad\forall
k\ge \bar k,\ \forall 
D\in\cC_{1+\bar \e}(\overline\Omega,\overline\Omega).
\eeq
Moreover, from Proposition \ref{Pn2.4} we infer that, if we choose
$\bar k$ large enough, for all $k\ge\bar k$ and $z\in Z_k$ the
function $u^{T^{\bar L}_k}_{k,z}$ is a solution of the Dirichlet problem
\beq
\label{2.25}
-\Delta u=|u|^{p-1}u+w\quad\mbox{ in }T_k^{\bar L}\left({1\over k}\,
  C_z\right),\qquad u=0\quad\mbox{ on }\partial  T_k^{\bar L}\left({1\over k}
  \,C_z\right).
\eeq
Thus, taking into account Proposition \ref{Pn2.5} we have to prove
that $E'(u_k)[v\cdot Du_k]=0$ for all vector field
$v\in\cC^1(\Omega,\R^n)$ such that $v\cdot\nu=0$ on $\partial\Omega$.

\no Therefore, for all vector field $v\in \cC^1(\overline\Omega, \R^n)$
such that $v\cdot\nu=0$ on $\partial \Omega$ and for all $\tau\in\R$,
let us consider the function
$D_\tau:\overline\Omega\to\overline\Omega$ defined by the Cauchy
problem
\beq
{\partial D_\tau(x)\over \partial \tau}=v\circ D_\tau(x),\qquad
D_0(x)=x\qquad\forall \tau\in\R,\ \forall x\in\overline\Omega.
\eeq
Then, we have $D_\tau(\overline\Omega)=\overline\Omega$
$\forall\tau\in\R$ and 
\beq
\lim_{\tau\to0}\cL(D_\tau\circ T^{\bar L}_k)=\cL(T^{\bar L}_k),
\eeq
so there exists $\bar\tau>0$ such that $D_\tau\circ T^{\bar L}_k\in
\cC_{\bar L}(P_k,\overline\Omega)$ $\forall
\tau\in[-\bar\tau,\bar\tau]$. 
It follows that 
\beq
\label{le}
E(u_{k})=\sum_{z\in Z_k}E(u^{T_k^{\bar L}}_{k,z})
\le \sum_{z\in Z_k}E(u^{D_\tau\circ T_k^{\bar L} }_{k,z})
= E(u^{D_\tau\circ T_k^{\bar L} }_{k})\quad\forall
\tau\in[-\bar\tau,\bar\tau].
\eeq
Notice that 
\beq
{d\phantom{1}\over d\tau}\, E(u_k\circ
D^{-1}_\tau)_{|_{\tau=0}}=-E'(u_k)[v\cdot \D u_k].
\eeq
Thus, we have to prove that
\beq
\label{=0}
{d\phantom{1}\over d\tau}\, E(u_k\circ D^{-1}_\tau)_{|_{\tau=0}}=0.
\eeq
For the proof, we argue by contradiction and assume that (\ref{=0})
does not hold.
For example, we assume that
\beq
\label{<0}
{d\phantom{1}\over d\tau}\, E(u_k\circ D^{-1}_\tau)_{|_{\tau=0}}<0
\eeq
(otherwise we replace $v$ by $-v$).
As a consequence, there exists a sequence of positive numbers
$(\tau_i)_{i\in\N}$ such that $\lim\limits_{i\to\infty}\tau_i=0$ and
$E(u_k\circ D^{-1}_{\tau_i})<E(u_k)$ $\forall i\in\N$.
From Corollary \ref{C2.2} we infer that, if we choose $\bar k$ large
enough, for all $k\ge\bar k$, $z\in Z_k$ and $i\in\N$ there exists a
unique minimizing function $\tilde u^{D_{\tau_i}\circ T^{\bar
    L}_k}_{k,z}$ and $\tilde u^{D_{\tau_i}\circ T^{\bar
    L}_k}_{k,z}\to \tilde u_{k,z}^{T^{\bar L}_k}$ as $i\to\infty$ in
$H^1_0(\Omega)$ $\forall k\ge \bar k$, $\forall z\in Z_k$.

\no As in the proof of Proposition \ref{Pn2.4}, let us consider the
functions $G^i_{k,z}$ verifying
\beq
\begin{array}{ll}G^i_{k,z}(x,t) =
{|t|^{p+1}\over p+1}+w(x)\, t &{\rm if }\ \sigma(z)[t-\tilde u^{D_{\tau_i}\circ T^{\bar
    L}_k}_{k,z}(x)]\ge 0,
\vspace{3mm}
\\
{\partial^2  G^i_{k,z}\over \partial t^2}\, (x,t)={\partial^2  G^i_{k,z}\over \partial t^2}\,
(x,\tilde u^{D_{\tau_i}\circ T^{\bar
    L}_k}_{k,z}(x)) &{\rm otherwise},
\end{array}
\eeq
the functional $E^i_{k,z}:H^1_0\left(T^{\bar L}_k\left({1\over k}\,
    C_z\right)\right)\to\R$ defined by
\beq
E_{k,z,i}(u)={1\over 2}\int_{T^{\bar L}_k\left({1\over k}\,
    C_z\right)}|\D u|^2dx-\int_{T^{\bar L}_k\left({1\over k}\,
    C_z\right)}G^i_{k,z}(x,u)\, dx
\eeq
and the manifold
\beq
\Gamma^i_{k,z} =\left\{u\in  H^1_0\left( T^{\bar L}_k\left({1\over k}\, C_z\right)\right)
\ :\ 
u\not\equiv \tilde u^{D_{\tau_i}\circ T^{\bar
    L}_k}_{k,z},\ E_{k,z,i}'(u)[u-\tilde u^{D_{\tau_i}\circ T^{\bar
    L}_k}_{k,z}]=0\right\}.
\eeq
We say that 
\beq
\label{max<}
\max\left\{\sum_{z\in Z_k}E\left(\tilde u^{T^{\bar L}_k}_{k,z}\circ
D^{-1}_{\tau_i}+t_z \left( u^{T^{\bar L}_k}_{k,z}\circ
D^{-1}_{\tau_i}-\tilde u^{T^{\bar L}_k}_{k,z}\circ
D^{-1}_{\tau_i}\right)\right) : t_z\ge 0,\  \forall z\in Z_k\right\}<E(u_k)
\eeq
for $i$ large enough.

\no In fact, arguing by contradiction, assume that (up to a
subsequence still denoted by $(\tau_i)_{i\in\N}$) the inequality
(\ref{max<}) does not hold.
Then, for all $i\in\N$ and $z\in Z_k$, there exists $t_{z,i}\ge 0$
such that 
\beq
\sum_{z\in Z_k}E\left(\tilde u^{T^{\bar L}_k}_{k,z}\circ
D^{-1}_{\tau_i}+t_{z,i} \left(u^{T^{\bar L}_k}_{k,z}\circ
D^{-1}_{\tau_i}- \tilde u^{T^{\bar L}_k}_{k,z}\circ
D^{-1}_{\tau_i}\right)\right)\ge E(u_k)\qquad\forall i\in\N.
\eeq
It follows that $\lim\limits_{i\to\infty}t_{z,i}=1$ $\forall z\in Z_k$ and
$$
\hspace{-3cm}\sum_{z\in Z_k}E\left(\tilde u^{T^{\bar L}_k}_{k,z}\circ
D^{-1}_{\tau_i}+t_{z,i} \left(u^{T^{\bar L}_k}_{k,z}\circ
D^{-1}_{\tau_i}- \tilde u^{T^{\bar L}_k}_{k,z}\circ
D^{-1}_{\tau_i}\right)\right)
$$
\beq
\hspace{3cm}\ge 
 \sum_{z\in Z_k}E\left(\tilde u^{T^{\bar L}_k}_{k,z}+
t_{z,i} \left(u^{T^{\bar L}_k}_{k,z}- \tilde u^{T^{\bar L}_k}_{k,z}\right)\right)
\qquad\forall i\in \N
\eeq
which, as $i\to\infty$, implies
\beq
\label{ge0}
{d\phantom{1}\over d\tau}\, E\left(u^{T^{\bar L}_k}_{k,z}\circ
D^{-1}_\tau\right)_{|_{\tau=0}}\ge 0 
\eeq
in contradiction with (\ref{<0}).
Thus, (\ref{max<}) holds.

\no Notice that, if $\bar k$ is chosen large enough, 
$$
\hspace{-3cm}E'\left(\tilde u^{T^{\bar L}_k}_{k,z}\circ
D^{-1}_{\tau_i}+t \left(u^{T^{\bar L}_k}_{k,z}\circ
D^{-1}_{\tau_i}- \tilde u^{T^{\bar L}_k}_{k,z}\circ
D^{-1}_{\tau_i}\right)\right)\cdot
$$
\beq
\hspace{3cm}\cdot \left[
\tilde u^{T^{\bar L}_k}_{k,z}\circ D^{-1}_{\tau_i}
-\tilde u^{D_{\tau_i}\circ T^{\bar L}_k}_{k,z}
+t \left(u^{T^{\bar L}_k}_{k,z}\circ D^{-1}_{\tau_i}-
\tilde u^{T^{\bar L}_k}_{k,z}\circ D^{-1}_{\tau_i}\right)\right],
\eeq
for $i$ large enough, is positive for  $t=\|u^{T^{\bar
    L}_k}_{k,z}\circ D^{-1}_{\tau_i}- \tilde u^{T^{\bar
    L}_k}_{k,z}\circ D^{-1}_{\tau_i}\|^{-1}_{H^1_0}$ and tends to
$-\infty$ as $t\to\infty$.
As a consequence, there exists $t^i_{k,z}\in\R$ such that 
\beq
\label{2.37}
\tilde u^{T^{\bar L}_k}_{k,z}\circ
D^{-1}_{\tau_i}+t^i_{k,z} \left(u^{T^{\bar L}_k}_{k,z}\circ
D^{-1}_{\tau_i}- \tilde u^{T^{\bar L}_k}_{k,z}\circ
D^{-1}_{\tau_i}\right)\in\Gamma^i_{k,z}.
\eeq
Therefore, from (\ref{max<}) and (\ref{2.37}) we obtain
\begin{eqnarray}
\nonumber 
E(u_k)
& > &\hspace{-3mm}
\max \left\{ \sum_{z\in Z_k}E\left(\tilde u^{T^{\bar L}_k}_{k,z}\circ
D^{-1}_{\tau_i}+t_z \left(u^{T^{\bar L}_k}_{k,z}\circ
D^{-1}_{\tau_i}- \tilde u^{T^{\bar L}_k}_{k,z}\circ
D^{-1}_{\tau_i}\right)\right) : t_z\ge 0\ \forall z\in Z_k\right\}
\\
&\ge&\hspace{-3mm}
 \sum_{z\in Z_k}E\left( u^{D_{\tau_i}\circ T_k^{\bar
      L}}_{k,z}\right)= E(u^{D_{\tau_i}\circ T_k^{\bar L}}_{k} )
\end{eqnarray}
for $i$ large enough, in contradiction with (\ref{le}).

\no Thus, we can conclude that
$
{d\phantom{1}\over d\tau}E(u_k\circ D^{-1}_\tau)_{|_{\tau=0}}=0
$
that is $E'(u_k)[v\cdot Du_k]=0$ for all vector field $v\in
\cC^1(\overline\Omega,\R^n)$ such that $v\cdot\nu=0$ on $\partial
\Omega$, so $u_k$ is a solution of problem (\ref{P}).

\no Notice that, if $J(k)$ denotes the number of elements of $Z_k$,
the solution $u_k$ has at least $J(k)$ nodal regions for $k$ large
enough.
Moreover, we have 
\beq
\lim_{k\to\infty} {J(k)\over k^n}=\meas (\Omega),
\eeq
so the number of nodal regions of $u_k$ tends to infinity as
$k\to\infty$.

\no Finally, notice that (\ref{ktoinfty}) follows directly from Lemma
\ref{Ln2.2} and Lemma \ref{Ln2.3}.

\no So the proof is complete.

\qed

\no Let us point out that, if $n=1$, condition (\ref{2.20}) in Theorem
\ref{MT} is satisfied.
In fact, it is a consequence of the following lemma
(see also Remark \ref{Rn3.1} concerning the case $n>1$).

\begin{lemma}
\label{ML}
Assume $n=1$, $p>1$ and $w\in L^2(\Omega)$.
Then, for all $L>1$ there exists $\bar k(L)\in \N$ such that 
\beq
\label{a}
\{T\in\cC_L(P_k,\overline\Omega)\ :\
T(P_k)=\overline\Omega\}\neq\emptyset\qquad\forall k\ge\bar k(L).
\eeq
Moreover,
\beq
\label{b}
\lim_{k\to\infty} \cL(T^L_k)=1\qquad\forall L>1.
\eeq
\end{lemma}

\no\proof Let $\Omega=]a,b[$. 
First notice that, if $Z_k$ consists of $J(k)$ points
$z_1,\ldots,z_{J(k)}$, then ${J(k)\over k}\le b-a$ and
$\lim\limits_{k\to\infty}{J(k)\over k}=b-a$.
Since $L>1$, it follows that there exists $\bar k(L)\in\N$ such that
${J(k)\over k}>{b-a\over L}$ $\forall k\ge \bar k(L)$, which implies
(\ref{a}) as one can easily verify.
Also, notice that in this case we have 
\beq
\label{*}
\hspace{-3cm} \cL(T^L_k)=\min\left\{\cL\ :\ \cL\ge 1,\ 
{1\over \cL}\le k\,
  \meas\left[T^L_k\left({1\over k}\,
      C_{z_j}\right)\right]\le\cL\right.
\eeq
$$
\left. \phantom{\int_\Omega}\hspace{7cm}  \mbox{ for } j=1,\ldots,J(k)\right\}\qquad \forall k > \bar k(L).
$$
Moreover, if we denote by $\bar u^{T^L_k}_{k,z}$ the function
$u^{T^L_k}_{k,z}$ obtained when $w=0$, we get by direct computation
\beq
\limsup_{k\to\infty} {1\over k^{p+3\over p-1}}\, \max\left\{E(\bar
  u^{T^L_k}_{k,z})\ :\ z\in Z_k\right\}<\infty
\eeq
and
\beq
\liminf_{k\to\infty} {1\over k^{p+3\over p-1}}\, \min\left\{E(\bar
  u^{T^L_k}_{k,z})\ :\ z\in Z_k\right\}>0.
\eeq
Taking into account that 
\beq
\hspace{-5.5cm}
\left|\int_{T^L_k\left({1\over k}\, C_z\right)}u\, w\, dx\right|
\le
\left(\int_{T^L_k\left({1\over k}\, C_z\right)}u^2\, dx\right)^{1\over
  2}\|w\|_{L^2(\Omega)}
\eeq
$$
\hspace{3.7cm}\le
\left[\meas T^L_k\left({1\over k}\, C_{z}\right)\right]^{{1\over
    2}-{1\over p+1}}\|u\|_{L^{p+1}\left(T^L_k\left({1\over k}\,
      C_z\right)\right)}\|w\|_{L^2(\Omega)}\quad \forall z\in Z_k,
$$
for all $w\in L^2(\Omega)$ we obtain
\beq
\limsup_{k\to\infty}{1\over k^{5-p\over 2(p-1)}}\max\{|E( u^{T^L_k}_{k,z})-
E(\bar u^{T^L_k}_{k,z})|\ :\ z\in Z_k\}<\infty
\eeq
and, as a consequence, 
\beq
\label{r}
\lim_{k\to\infty}  {J(k)\over k^{p+3\over 2(p-1)}} \max \{|E( u^{T^L_k}_{k,z})-
E(\bar u^{T^L_k}_{k,z})|\ :\ z\in Z_k\}<\infty.
\eeq
It is clear that $1\le\cL(T^L_k)\le L$, so
$\liminf\limits_{k\to\infty}\cL(T^L_k)\ge 1$.
Arguing by contradiction, assume that (\ref{b}) does not hold, that is
\beq
\limsup_{k\to\infty}\cL(T^L_{k})>1.
\eeq
Thus, there exists a sequence $(k_i)_{i\in\N}$ in $\N$ such that 
\beq
\lim_{i\to\infty}\cL(T^L_{k_i})>1.
\eeq
Taking into account (\ref{*}), it follows that there exists a sequence
$(z'_i)_{i\in\N}$ such that $z'_i\in Z_{k_i}$ $\forall i\in \N$ and
(up to a subsequence)
\beq
\label{A}
\lim_{i\to\infty}k_i\meas\left[T^L_{k_i}\left({1\over k_i}\,
    C_{z'_i}\right)\right]>1,
\eeq
or there exists a sequence $(z_i'')_{i\in\N}$ such that $z_i''\in
Z_{k_i}$ $\forall i\in \N$ and (up to a subsequence)
\beq
\label{B}
\lim_{i\to\infty}k_i\meas\left[T^L_{k_i}\left({1\over k_i}\,
    C_{z_i''}\right)\right]<1.
\eeq
Let us consider first the case where (\ref{A}) holds.
In this case, for all $i\in\N$, choose $\tilde z_i'$ in $Z_{k_i}$ such
that 
\beq
\meas\left[T^L_{k_i}\left({1\over k_i}\, C_{\tilde z_i'}\right)\right]
=
\min\left\{
\meas\left[T^L_{k_i}\left({1\over k_i}\, C_{z}\right)\right]\ :
\ z\in Z_{k_i}\right\}.
\eeq
Then, taking into account that 
\beq
(b-a)=\sum_{z\in Z_{k_i}}\meas\left[T^L_{k_i}\left({1\over k_i}\,
    C_{z}\right)\right]\ge J(k_i)\meas\left[T^L_{k_i}\left({1\over
      k_i}\, C_{\tilde z_i'}\right)\right]
\eeq
and that $\lim\limits_{i\to\infty}{J(k_i)\over k_i}=b-a$, we obtain
\beq
\label{le1}
\limsup_{i\to\infty}k_i\cdot \meas\left[T^L_{k_i}\left({1\over
      k_i}\, C_{\tilde z_i'}\right)\right]\le 1.
\eeq
As one can easily verify, for all $i\in\N$, there exists a function
$T'_i\in\cC_L(P_{k_i},\overline\Omega)$ such that
$T'_i(P_{k_i})=\overline\Omega$, 
\beq
\label{1}
\meas\left[T'_i\left({1\over
      k_i}\, C_{ z}\right)\right]=\meas\left[T^L_{k_i}\left({1\over
      k_i}\, C_{z}\right)\right]\qquad\forall z\in
Z_{k_i}\setminus\{z'_i,\tilde z'_i\},
\eeq
\beq
\label{2}
\meas\left[T'_i\left({1\over
      k_i}\, C_{z_i'}\right)\right]
=
\meas\left[T'_i\left({1\over
      k_i}\, C_{\tilde z_i'}\right)\right]
\eeq
and
\beq
\label{3}
\lim_{i\to\infty}\max\{|T'_i(x)-T^L_{k_i}(x)|\ :\ x\in P_{k_i}\}=0.
\eeq
Taking into account that $E(\bar u^{T_i'}_{k_i,z})=E(\bar
u^{T_{k_i}^L}_{k_i,z})$ $\forall z\in Z_{k_i}\setminus\{z_i',\tilde
z_i'\}$, from (\ref{r}) we infer that
\beq
\lim_{i\to\infty}{1\over k_i^{p+3\over p-1}}\sum_{z\in
  Z_{k_i}\setminus\{ z_i',\tilde z_i'\}}|E(u^{T_i'}_{k_i,z}) -E(
u^{T_{k_i}^L}_{k_i,z})
|=0.
\eeq
On the other hand,
\beq
\liminf_{i\to\infty}{1\over k_i^{p+3\over p-1}}[E(u^{T_i'}_{k_i,z_i'}) +E(
u^{T_i'}_{k_i,\tilde z_i'})]>0
\eeq
as one can verify by direct computation.
Moreover, (\ref{A}) and (\ref{le1}) imply
\beq
\label{>}
\liminf_{i\to\infty}
{E(u^{T^L_{k_i}}_{k_i,z_i'}) +E( u^{T^L_{k_i}}_{k_i,\tilde z_i'})
\over
E(u^{T_i'}_{k_i,z_i'}) +E( u^{T_i'}_{k_i,\tilde z_i'})}
>1.
\eeq
It follows that 
\beq
\liminf_{i\to\infty}{1\over k_i^{p+3\over p-1}}
\sum_{z\in Z_{k_i}}
[E(u^{T^L_{k_i}}_{k_i,z}) -E( u^{T'_{i}}_{k_i, z})]>0,
\eeq
which is a contradiction because
\beq
\sum_{z\in Z_{k_i}} E(u^{T^L_{k_i}}_{k_i,z}) 
\le
\sum_{z\in Z_{k_i}} E( u^{T'_{i}}_{k_i, z})
\qquad\forall i\in\N.
\eeq
When the case (\ref{B}) occurs, we argue in analogous way.
In this case, for all $i\in\N$ we choose $\tilde z'' $ in $Z_{k_i}$
such that 
\beq
\meas\left[T^L_{k_i}\left({1\over k_i}\, C_{\tilde z_i''}\right)\right]
=
\max\left\{
\meas\left[T^L_{k_i}\left({1\over k_i}\, C_{z}\right)\right]\ :
\ z\in Z_{k_i}\right\},
\eeq
which implies
\beq
\label{ge1}
\liminf_{i\to\infty}k_i\cdot \meas\left[T^L_{k_i}\left({1\over
      k_i}\, C_{\tilde z_i''}\right)\right]\ge 1.
\eeq
Moreover, we can consider a function
$T''_i\in\cC_L(P_{k_i},\overline\Omega)$ satisfying all the properties
of $T'_i$ with $z_i''$ and $\tilde z_i''$ instead of $z_i'$ and
$\tilde z_i'$.

\no Then, we can repeat for $T_i''$, $z_i''$ and $\tilde z_i''$ the
same arguments as before.
In particular, the property 
\beq
\liminf_{i\to\infty}
{E(u^{T^L_{k_i}}_{k_i,z_i''}) +E( u^{T^L_{k_i}}_{k_i,\tilde z_i''})
\over
E(u^{T_i''}_{k_i,z_i''}) +E( u^{T_i''}_{k_i,\tilde z_i''})}
>1
\eeq
(analogous to (\ref{>})) now follows from (\ref{B}) and (\ref{ge1}).
Thus, also in this case we obtain again a contradiction with the
minimality property of $\sum\limits_{z\in
  Z_{k_i}}E(u^{T^L_{k_i}}_{k_i,z})$.
So the proof is complete.

\qed

\no As a direct consequence of Theorem \ref{MT}, Proposition \ref{MP}
and Lemma \ref{ML} we obtain the following corollary.
\begin{cor}
Assume $n=1$ and $p>1$.
Then, for all $w\in L^2(\Omega)$, problem (\ref{P}) has infinitely
many solutions.
More precisely, for all $L>1$ there exists $\bar k(L)\ge k_\Omega$
such that for all $k\ge \bar k(L)$ there exists $T^L_k\in
\cC_L(P_k,\overline\Omega)$ such that $T^L_k(P_k)=\overline\Omega$ and
the function $u_k=\sum\limits_{z\in Z_k} u^{T^L_k}_{k,z}$ is a
solution of problem (\ref{P}).

\no Moreover, the number of nodal regions of $u_k$ tends to infinity as
$k\to\infty$ and
\beq
\lim_{k\to\infty}\cL(T^L_k)=1,\quad
\lim_{k\to\infty}\min\{E(u^{T^L_k}_{k,z})\ :\ z\in
Z_k\}=\infty\quad\forall L>1.
\eeq
\end{cor}


\sezione{Final remarks}



Notice that the method we used in Section \ref{S2} to find infinitely many
solutions of problem (\ref{P}) with a large number of nodal regions
having a prescribed structure (a check structure) may be used also in
other elliptic problems as we show in this section.

\no It is clear that in this method condition (\ref{2.20}) plays a
crucial role.
In Section \ref{S2}  this condition is proved only in the case $n=1$.
In next remark, we discuss about the case $n>1$.
 
\begin{rem}
\label{Rn3.1}
{\em
Assume that condition (\ref{2.20}) does not hold.
Then, there exists a sequence $(L_i)_{i\in\N}$ in $\R$ such that
\beq
\lim\limits_{i\to\infty}L_i=\infty\qquad{\rm and }\qquad
\limsup_{k\to\infty}\cL(T^{L_i}_{k})=L_i\quad\forall i\in\N.
\eeq
As a consequence, we can construct a sequence $(k_i)_{i\in\N}$ such
that
\beq
\label{itoinfty}
\lim_{i\to\infty}k_i=\infty,\qquad\lim_{i\to\infty}\min\big\{E\big(u^{T^{L_i}_{k_i}}_{k_i,z}\big)\
:\ z\in Z_{k_i}\big\}=\infty,\qquad
\lim_{i\to\infty}\cL(T^{L_i}_{k_i})=\infty.
\eeq
Notice that $\cL(T^{L_i}_{k_i})$ is large, for example, when there are large
differences in the sizes or in the shapes of the subdomains
$T^{L_i}_{k_i}\left({1\over k_i}\, C_z\right)$ with $z\in Z_{k_i}$.
For $k_i$ large enough, too large differences seem to be incompatible
with the minimality property
\beq
\sum_{z\in Z_{k_i}}E\big(u^{T^{L_i}_{k_i}}_{k_i,z}\big)=\min\left\{
\sum_{z\in Z_{k_i}}E\big(u^T_{k_i,z}\big)\ :\
  T\in\cC_{L_i}(P_{k_i},\overline\Omega),\
  T(P_{k_i})=\overline\Omega\right\}\quad\forall i\in\N.
\eeq
This fact explains why condition \eqref{2.20} holds in the case $n=1$.
In the case $n>1$, on the contrary, even if the subdomains $T^{L_i}_{k_i}\left(\frac{1}{k_i}C_z\right)$ with $z\in Z_{k_i}$ have all the same shape and the same size, we cannot exclude that $\cL\left(T^{L_i}_{k_i}\right)$ is large as a consequence of the fact that the shape of these subdomains is very different from the cubes of $\R^n$.
This explains why it is difficult to prove that condition \eqref{2.20} holds also for $n>1$.

Therefore, in the case $n>1$, the natural idea is to restrict the class of the 
admissible deformations of the nodal regions.

For example, we can fix $\tilde L\ge L'_\Omega$, $T_0\in \cC_{\tilde L}(\overline\Omega,\overline\Omega)$, $r>0$ and consider the set of deformations
\beq
\label{R1}
\Theta^{\tilde L}_k(T_0,r)=\{T\in \cC_{\tilde L}(P_k,\overline\Omega)\ :\ T(P_k)=\overline\Omega,\ d_k(T,T_0)\le r\}
\eeq
where 
\beq
d_k(T,T_0)=\sup_{P_k}|T-T_0|+\Lip (T-T_0)
\eeq
with
\beq
\Lip (T-T_0)=\sup\{|x-y|^{-1}|T(x)-T_0(x)-T(y)+T_0(y)|\ :\ x,y \text{ in }P_k,\ x\neq y\}.
\eeq
Then, arguing exactly as in Section 2 but minimizing in the subset $\Theta^{\tilde L}_k(T_0,r)$ (instead of the set \eqref{Fdiv0}), we obtain a minimizing deformation $T^{\tilde L,r}_k$ which, for $k$ large enough, gives rise to a solution $u^{\tilde L,r}_k$ of problem \eqref{P} provided the condition 
\beq
\label{R*}
\limsup_{k\to\infty}\cL(T^{\tilde L,r}_k)<\tilde L,
\quad 
\limsup_{k\to\infty}d_k(T^{\tilde L,r}_k, T_0)<r
\eeq
(analogous to condition \eqref{2.20}) is satisfied.

It is clear that condition \eqref{R*} holds or fails depending on the choice of $\tilde L,T_0$ and $r$ that have to be chosen in a suitable way.
For example, in the case $n=1$, if we choose $T_0(x)=x$ $\forall x\in\Omega$, \eqref{R*} holds for all $\tilde L>1$ and $r>0$ as follows from Lemma \ref{ML}.

In the case $n>1$, condition \eqref{R*} seems to have more chances than condition \eqref{2.20} to be satisfied.
In fact, as we show in a paper in preparation, a variant of this method works for example when $\Omega$ is a cube of $\R^n$ with $n>1$, $p>1$, $p<\frac{n+2}{n-2}$ if $n>2$ and, for all $w\in L^2(\Omega)$, allows us to find infinitely many solutions $u_k(x)$ such that the nodal regions of $u_k\left(\frac{x}{k}\right)$, after translations, tend to the cube as $k\to\infty$.

Therefore, it seems quite natural to expect that, by a suitable choice of $\tilde L$, $T_0$ and $r$, for every bounded domain $\Omega$ in $\R^n$ with $n>1$ and for all $w\in L^2(\Omega)$ one can find infinitely many nodal solutions of problem \eqref{P} with $p>1$ and $p<\frac{n+2}{n-2}$ if $n>2$.

\hfill $\Box$
}
\end{rem}

\no Notice that this method to construct solutions with nodal regions
having this check structure works for more general nonlinearities,
even when they are not perturbations of symmetric nonlinearities: for
example when in problem (\ref{P}) the term $|u|^{p-1}u+w$ is replaced
by $c_+(u^+)^p-c_-(u^-)^p+w$ with $c_+>0$ and $c_->0$.

\no In fact, this method does not require any technique of deformation
from the symmetry.
For example, let us show how Lemma \ref{ML} has to be modified in this
case.

\no In this case the energy functional is 
\beq
F(u)={1\over 2}\int_\Omega|\D u|^2dx-{c_+\over
  p+1}\int_\Omega(u^+)^{p+1}dx-{c_-\over
  p+1}\int_\Omega(u^-)^{p+1}dx-\int_\Omega wu\, dx.
\eeq
We denote by $F_0$ the functional $F$ when $w=0$.

\no Now, consider the number $\hat L\ge 1$ defined by $\hat L={1\over
  \min\{\hat t,2-\hat t\}}$ where $\hat t\in]0,2[$ is the unique
number such that 
$$
\hspace{-4cm} \min\{F_0(u)\ :\ u\in H^1_0(]0,2[),\ u\ge 0\ {\rm in }\ ]0,\hat t[,\ u\le 0\ {\rm in }\
]\hat t,2[,
$$
\beq
 u^+\not\equiv 0,\ u^-\not\equiv 0,\ F'_0(u)[u^+]=0,\
F'_0(u)[u^-]=0\}
\eeq
$$
=\min\{F_0(u)\ :\ u\in H^1_0(]0,2[),\  u^+\not\equiv 0,\ u^-\not\equiv 0,\ F'_0(u)[u^+]=0,\
F'_0(u)[u^-]=0\}.
$$
Notice that $\hat t=1$ (and so $\hat L=1$) if and only if $c_+=c_-$.

\no Then, we have the following lemma which extends Lemma \ref{ML}.

\begin{lemma}
\label{Ln3.3}
Assume $n=1$, $c_+>0$, $c_->0$, $p>1$.
Then, for all $L>\hat L$ there exists $\hat k(L)\in\N$ such that
$\{T\in\cC_L(P_k,\overline\Omega)\ :\ T(P_k)=\overline\Omega\}\neq\emptyset$ $\forall k\ge\hat
k(L)$ and $\lim\limits_{k\to\infty}\cL(T^L_k)=\hat L$.
\end{lemma}

\no \proof Here we describe only how the proof of Lemma \ref{ML} has
to be modified in order to be adapted in this case.

\no First notice that, since $L>\hat L$ and $\hat L\ge 1$, there
exists $\hat k(L)\in\N$ such that 
\beq
\{T\in\cC_L(P_k,\overline\Omega)\ :\
T(P_k)=\overline\Omega\}\neq\emptyset\qquad\forall k\ge\hat k(L)
\eeq
as one can verify as in the proof of Lemma \ref{ML}.

\no In order to prove that $\lim\limits_{k\to\infty}\cL(T^L_k)=\hat
L$, we argue by contradiction and assume that there exists a sequence
$(k_i)_{i\in\N}$ in $\N$ such that $\lim\limits_{i\to\infty} k_i=\infty$ and
$\lim\limits_{i\to\infty}\cL(T^L_{k_i})\neq\hat L$.

\no First, notice that the case
\beq
\label{d<}
\lim_{i\to\infty}\cL(T^L_{k_i})<\hat L
\eeq
cannot happen.
In fact, for all $i\in\N$ we can choose $\hat z_i$ and $\hat z_i+1$
in $Z_{k_i}$ such that (up to a subsequence)
\beq
\lim_{i\to\infty} k_i\left[\meas T^L_{k_i}\left({1\over k_i}\, C_{\hat
      z_i}\right)+\meas T^L_{k_i}\left({1\over k_i}\, C_{\hat
      z_i+1}\right)\right]\le 2.
\eeq
Taking into account the minimality property
\beq
\label{dmin}
\sum_{z\in Z_{k_i}}F(u^{T^L_{k_i}}_{k_i,z})=\min\left\{
\sum_{z\in Z_{k_i}}F(u^T_{k_i,z})\ :\ T\in\cC_L(P_k,\overline\Omega),\
T(P_{k_i})=\overline\Omega\right\}
\qquad\forall i\in \N,
\eeq
it follows that 
\beq
\lim_{i\to\infty} \frac{2\min\{\meas T^L_{k_i}\left({1\over k_i}\, C_{\hat
      z_i}\right),\, \meas T^L_{k_i}\left({1\over k_i}\, C_{\hat
      z_i+1}\right)\}}{ \meas T^L_{k_i}\left({1\over k_i}\, C_{\hat
      z_i}\right)+\meas T^L_{k_i}\left({1\over k_i}\, C_{\hat
      z_i+1}\right)}=\min\{\hat t,2-\hat t\}.
\eeq
As a consequence, we obtain
\beq
\lim_{i\to\infty} k_i\min\left\{\meas T^L_{k_i}\left({1\over k_i}\, C_{\hat
      z_i}\right),\, \meas T^L_{k_i}\left({1\over k_i}\, C_{\hat
      z_i+1}\right)\right\}\le \min\{\hat t,2-\hat t\}
\eeq
which implies
\beq
\lim_{i\to\infty}\cL(T^L_{k_i})\ge {1\over \min\{\hat t,2-\hat
  t\}}=\hat L.
\eeq
In order to prove that $\lim_{i\to\infty}\cL(T^L_{k_i})=\hat
L$,  arguing by contradiction, assume that\linebreak 
$\lim_{i\to\infty}\cL(T^L_{k_i})$ $>\hat L$. 

\no As a consequence, since
\beq
\cL(T^L_{k_i})=\min\left\{\cL\ :\ \cL\ge 1, \ {1\over \cL}\le k_i\meas
  T^L_{k_i}\left({1\over k_i}\, C_z\right)\le\cL\ \ \forall z\in
  Z_{k_i}\right\}, 
\eeq
there exists a sequence $(z'_i)_{i\in\N}$ such that $z'_i\in Z_{k_i}$
$\forall i\in\N$ and 
\beq
\label{fA}
\lim_{i\to\infty} k_i\meas T^L_{k_i}\left({1\over k_i}\, C_{z'_i}\right)>\hat L,
\eeq
or there exists a sequence $(z''_i)_{i\in\N}$ such that $z''_i\in
Z_{k_i}$ $\forall i\in\N$ and 
\beq
\label{fB}
\lim_{i\to\infty} k_i\meas T^L_{k_i}\left({1\over k_i}\,
  C_{z''_i}\right)<{1\over \hat L}. 
\eeq

\no Assume, for example, that $\hat t\le 1$ (otherwise we argue in a
similar way but with $\hat t$ replaced by $2-\hat t$). 
Then,  $\hat L={1\over \hat t}$ and, if $\hat t=1$, Lemma \ref{ML} applies.
Thus, it remains to consider the case $\hat t\in]0,1[$.

\no Consider first the case where (\ref{fA}) holds. Notice that there
exists a sequence $(\zeta'_i)_{i\in\N}$ such that $\zeta'_i\in
Z_{k_i}$ and $|z_i'-\zeta'_i|=1$ $\forall i \in\N$.

\no Then, the minimality property (\ref{dmin}) implies
\beq
\label{f1}
\lim_{i\to\infty} \frac{2\max\{\meas T^L_{k_i}\left({1\over k_i}\, C_{
      z'_i}\right),\, \meas T^L_{k_i}\left({1\over k_i}\, C_{ 
      \zeta'_{i}}\right)\}}{ \meas T^L_{k_i}\left({1\over k_i}\, C_{ 
      z'_i}\right)+\meas T^L_{k_i}\left({1\over k_i}\, C_{ 
      \zeta'_{i}}\right)}= 2-\hat t
\eeq
and, arguing as in the proof of Lemma \ref{ML} but with $\meas
T^L_{k_i}\left({1\over k_i}\, C_{  
      z'_i}\right)+\meas T^L_{k_i}\left({1\over k_i}\, C_{ 
      \zeta'_{i}}\right)$ instead of $\meas T^L_{k_i}\left({1\over k_i}\, C_{ 
      z'_i}\right)$
also
\beq
\label{f2}
\lim_{i\to\infty} k_i\left[\meas T^L_{k_i}\left({1\over k_i}\, C_{ 
      z'_i}\right)+\meas T^L_{k_i}\left({1\over k_i}\, C_{ 
      \zeta'_{i}}\right)\right]= 2.
\eeq
As a consequence of (\ref{f1}) and (\ref{f2}), we obtain
\beq
\label{f3}
\lim_{i\to\infty} k_i\max\left\{\meas T^L_{k_i}\left({1\over k_i}\, C_{
      z'_i}\right),\, \meas T^L_{k_i}\left({1\over k_i}\, C_{ 
      \zeta'_{i}}\right)\right\}= 2-\hat t
\eeq
which is a contradiction because 
\beq
\lim_{i\to\infty} k_i\meas T^L_{k_i}\left({1\over k_i}\,
  C_{z'_i}\right)>\hat L={1\over \hat t} 
\eeq
with ${1\over \hat t}>2-\hat t$ for $\hat t\in ]0,1[$.

\no Thus, we can conclude that the case (\ref{fA}) cannot happen.

\no In a similar way we argue in order to obtain a contradiction in
the case (\ref{fB}). 
In fact, assume that (\ref{fB}) holds.
Notice that there exists a sequence $(\zeta''_i)_{i\in\N}$ such that
$\zeta''_i\in Z_{k_i}$ and $|z''_i-\zeta''_i|=1$ $\forall i\in\N$. 

\no As before, the minimality property (\ref{dmin}) implies that 
\beq 
\lim_{i\to\infty} \frac{2\min\{\meas T^L_{k_i}\left({1\over k_i}\, C_{
      z''_i}\right),\, \meas T^L_{k_i}\left({1\over k_i}\, C_{ 
      \zeta''_{i}}\right)\}}{ \meas T^L_{k_i}\left({1\over k_i}\, C_{ 
      z''_i}\right)+\meas T^L_{k_i}\left({1\over k_i}\, C_{ 
      \zeta''_{i}}\right)}= \hat t
\eeq
and
\beq
\lim_{i\to\infty} k_i\left[\meas T^L_{k_i}\left({1\over k_i}\, C_{ 
      z''_i}\right)+\meas T^L_{k_i}\left({1\over k_i}\, C_{ 
      \zeta''_{i}}\right)\right]= 2.
\eeq
As a consequence, we infer that
\beq 
\lim_{i\to\infty} k_i\min\left\{\meas T^L_{k_i}\left({1\over k_i}\, C_{
      z''_i}\right),\, \meas T^L_{k_i}\left({1\over k_i}\, C_{ 
      \zeta''_{i}}\right)\right\}= \hat t,
\eeq
which is in contradiction with (\ref{fB}) because 
\beq
\lim_{i\to\infty} k_i\meas T^L_{k_i}\left({1\over k_i}\, C_{z''_i}\right)<{1\over \hat L}=\hat t.
\eeq
Thus, we can conclude that $\lim\limits_{i\to\infty}\cL(T^L_{k_i})=\hat L$ so the proof is complete.

\qed

\begin{rem}
{\em
\label{R0}
The results we present in this paper concern the existence of solutions with a large number of nodal regions.
In particular, when $\Omega\subset\R^n$ with $n=1$, these solutions must have, as a consequence, a large number of zeroes.
In next propositions we show that the term $w$ can be chosen in such a way that the sign of the solutions is related to the nodal regions of the eigenfunctions of the Laplace operator  $-\Delta$ in $H^1_0(\Omega)$.
In particular, if $n=1$ we show that for suitable terms $w$ in $L^2(\Omega)$, problem \eqref{P}
does not have solutions with a small number of zeroes: more precisely, we show that for all 
positive integer $h$ there exists $w_h\in L^2(\Omega)$ such that every solution of problem 
\eqref{P} has at least $h$ zeroes (it follows from Corollary \ref{Cr3}).
}
\end{rem}

\begin{lemma}
\label{Lr1}
Let $D$ be a piecewise smooth, bounded open subset of $\R^n$, $n\ge 1$, and $\lambda_1$ be the first eigenvalue of the Laplace operator $-\Delta$ in $H^1_0(D)$.

Let $g:D\times \R\to\R$ be a Carath\'eodory function such that 
\beq
\label{r1}
\inf\{g(x,t)-\lambda_1 t\ :\ x\in D, \ t\ge 0\}>0.
\eeq
Let $u\in H^1(D)$ be a weak solution of the equation
\beq
\label{r2}
-\Delta u=g(x,u)\quad\mbox{ in }D.
\eeq
Then $\inf_D u<0$. Moreover, if
\beq
\label{r3}
\sup\{g(x,t)-\lambda_1 t\ :\ x\in D, \ t\le 0\}<0,
\eeq
then $\sup_D u>0$. 
\end{lemma}

\proof Let $e_1$ be a positive eigenfunction corresponding to the eigenvalue $\lambda_1$, that is
\beq
\Delta e_1+\lambda_1 e_1=0,\quad e_1>0\ \mbox{ in } D,\quad e_1\in H^1_0(D).
\eeq
Arguing by contradiction, assume that \eqref{r1} holds and $u\ge 0$ in $D$.
Then, from \eqref{r2} we infer that
\beq
-\int_D\Delta u\, e_1\, dx=\int_D g(x,u)  e_1\, dx
\eeq
which implies
\beq
\int_D Du\, De_1\, dx=\int_{\partial D} u\, De_1\cdot\nu\, d\sigma -\int_D u\, \Delta e_1\, dx
=\int_D g(x,u) e_1\, dx,
\eeq
where $\nu$ denotes the outward normal on $\partial D$, so that
\beq
\int_{\partial D}u\, De_1\cdot\nu\, d\sigma\le 0,
\eeq 
and $g(x,t)\ge \lambda_1t+c$ $\forall x\in D$, $\forall t\in \R$ for a suitable constant $c>0$.
It follows that
\beq
\lambda_1\int_D u  e_1\, dx\ge\int_D g(x,u)  e_1dx\ge \lambda_1\int_Du  e_1\, dx+c\int_De_1\, dx,
\eeq
which implies $c\int_D e_1\, dx\le 0$, that is a contradiction.
Thus, the function $u$ cannot be a.e. nonnegative in $D$.

In a similar way one can show that we cannot have $u\le 0$ a.e. in $D$ when \eqref{r3} holds, 
so the proof is complete.

\qed

In particular, Lemma \ref{Lr1} may be used to obtain informations on the effect of the term $w$ on the sign changes of the solutions of problem \eqref{P}, as we describe in the following proposition.

\begin{prop}
\label{Pr2}
Let $\Omega\subset\R^n$ with $n\ge 1$ and $e_k\in H^1_0(\Omega)$ be an eigenfunction of the Laplace operator $-\Delta$ with eigenvalue $\lambda_k$, that is $\Delta e_k+\lambda_k e_k=0$ in $\Omega$.
Assume that $w\in L^2(\Omega)$ satisfies
\beq
\label{r4}
e_kw\ge 0\ \mbox{ in }\Omega,\quad\inf_\Omega |w|>\max\{\lambda_k t-t^p\ :\ t\ge 0\}.
\eeq
Let $\Omega _k\subseteq\Omega$ be a nodal region of $e_k$, that is $e_{k|_{\Omega_k}}\in H^1_0(\Omega_k)$, $e_k\neq 0$ everywhere in $\Omega_k$ and the sign of $e_k$ is constant in $\Omega_k$.

Then, there exists no function $u$ in $H^1(\Omega)$ satisfying 
\beq
u\, e_k\ge 0\ \mbox{ and }-\Delta u=|u|^{p-1}u+w\quad\mbox{ in }\Omega_k.
\eeq
In particular, if $w$ satisfies \eqref{r4}, every solution $u$ of problem \eqref{P} must satisfy 
\beq
\inf_{\Omega_k}u<0\ \mbox{  if }\ e_k>0\ \mbox{  in }\  \Omega_k\ \mbox{  and }\  \sup_{\Omega_k}u>0\ \mbox{  if }\  e_k<0\ \mbox{  in }\ \Omega_k.
\eeq
\end{prop}

\proof Notice that $\lambda_k$ is the first eigenvalue of the Laplace operator $-\Delta$ in $H^1_0(\Omega_k)$ and $|e_k|$ is a corresponding positive eigenfunction.
Moreover, if we set $g(x,t)=|t|^{p-1}t+w(x)$, we infer from \eqref{r4} that, if $w(x)>0$, 
\beq
g(x,t)\ge \lambda_kt+\tilde c\qquad\forall t\ge 0
\eeq
and, if $w(x)<0$,
\beq
g(x,t)\le \lambda_kt-\tilde c\qquad\forall t\le 0
\eeq
where
\beq
\tilde c=\inf_\Omega|w|-\max\{\lambda_k t-t^p\ :\ t\ge 0\}>0.
\eeq
Since $u e_k\ge 0$ and $w  e_k\ge 0$ in $\Omega_k$ and $e_k$ has constant sign in $\Omega_k$, we have $u\ge 0$ and $w>0$ in $\Omega_k$ if $e_k>0$ in $\Omega_k$ and $u\le 0$, $w<0$ in $\Omega_k$ in the opposite case. 
Thus, our assertion follows from Lemma \ref{Lr1}.
In fact, for example, in the case $e_k>0$ in $\Omega_k$ we cannot have
\beq
-\Delta u=|u|^{p-1}u+w\quad\mbox{ in }\Omega_k
\eeq
otherwise $\inf_{\Omega_k} u<0$, because of Lemma \ref{Lr1}, while $u\ge 0$ in $\Omega_k$.

In the opposite case, when $e_k<0$ in $\Omega_k$, one can argue in a similar way, so the proof is complete.

\qed

\begin{cor}
\label{Cr3}
Assume $n=1$ and $\Omega=]a,b[$.
Let us denote by $\lambda_1<\lambda_2<\lambda_3<\ldots$ the eigenvalues of the Laplace operator $-\Delta $ in $H^1_0(]a,b[)$ and, for all $k\in\N$, consider an eigenfunction $e_k$ with eigenvalue $\lambda_k$.
Moreover, assume that $w\in L^2(]a,b[)$ satisfy condition \eqref{r4}.

Set $h=k-1$ ($h$ is the number of zeroes of $e_k$ in $]a,b[$).

Then, every solution of problem \eqref{P} has in $]a,b[$ at least $h$ zeroes $\zeta_1,\zeta_2,\ldots,\zeta_h$ such that
\beq
\left|a+\frac{b-a}{k}\, i-\zeta_i\right|<\frac{b-a}{k}\qquad\mbox{ for }i=1,\ldots,h.
\eeq
\end{cor}

\proof Notice that the points $\nu_i=a+\frac{b-a}{k}\, i$, for $i=0,1,\ldots,k$, are the zeroes of $e_k$ in $[a,b]$ and the intervals 
$I_i=]\nu_{i-1},\nu_i[$, for $i=1,\ldots,k$, are the nodal regions of $e_k$.

Assume, for example, that $e_k>0$ on $I_1$ (in a similar way one can argue if $e_k<0$ in $I_1$).
Then, from Proposition \ref{Pr2} we infer that for every solution $u$ of problem \eqref{P} we have $\inf_{I_i}u<0$ for $i$ odd and 
$\sup_{I_i}>0$ for $i$ even.

Therefore, the function $u$ has at least $h$ zeroes $\zeta_1,\ldots,\zeta_h$ such that $\zeta_i\in ]\nu_{i-1},\nu_{i+1}[$ for 
$i=1,\ldots,h$, so the proof is complete.

\qed

\begin{rem}
\label{Rr4}
{\em
Notice that all the assertions in Proposition \ref{Pr2} and Corollary \ref{Cr3} still hold when the nonlinear term $|u|^{p-1}u$ is replaced by $c_+(u^+)^p-c_-(u^-)^p$ where $c_+$ and $c_-$ are two positive constants.
In this case we have only to replace $\max\{\lambda_kt-t^p\ :\ t\ge 0\}$ by $\max\{\lambda_k t-\bar c t^p\ :\ t\ge 0\}$, where $\bar c=\min\{c_+,c_-\}>0$.
}\end{rem}

\no Notice that this method to construct solutions with nodal regions
having a check structure may be used for nonlinear elliptic problems
with different boundary conditions, for systems and also when the
nonlinear term has critical growth. 
For example, for all $\lambda\in\R$ consider the Dirichlet problem 
\beq
\label{e*}
-\Delta u=|u|^{4\over n-2}u+\lambda u+w\quad\mbox{ in
}\Omega,\qquad u=0\quad\mbox{ on
}\partial\Omega
\eeq
whose solutions are critical points of the energy functional
$\cF:H^1_0(\Omega)\to \R$ defined by
\beq
\cF(u)={1\over 2}\int_\Omega|\D u|^2dx-{n-2\over 2n}\int_\Omega
|u|^{2n\over n-2}dx-{\lambda\over 2} \int_\Omega
u^2dx-\int_\Omega w\, u\, dx\qquad u\in H^1_0(\Omega).
\eeq
Using this method, if the functional $\cF$ satisfies condition
(\ref{2.20}), one can prove that for $n\ge 4$ and $\lambda>0$ the
functional $\cF$ has an unbounded sequence of critical levels.
More precisely, the following theorem can be proved.

\begin{teo}
\label{Tn3.4}
Let $n\ge 4$, $\lambda>0$, $w\in L^2(\Omega)$ and assume that
condition  (\ref{2.20}) holds for the functional $\cF$.
Then, there exists $\bar k\ge k_\Omega$ such that for all $k\ge\bar k$
there exists $T^{\bar L}_k\in\cC_{\bar L}(P_k,\overline\Omega)$ and a
solution $u_k$ of problem (\ref{e*}) such that $T^{\bar
  L}_k(P_k)=\overline\Omega$ and, if for all $z\in Z_k$ we set $u^z_k(x)=u_k(x)$ when
$x\in T^{\bar L}_k\left({1\over k}\, C_z\right)$, $u^z_k(x)=0$
otherwise, then we have $[\sigma(z)u^z_k]^+\not\equiv 0$,
\beq
 \cF(u^z_k)\le {1\over n}\, S^{n/2}\quad\forall z\in Z_k
\quad\mbox{ and }\quad \lim_{k\to\infty} \min\left\{
\cF(u^z_k)\ :\  z\in Z_k\right\}= {1\over n}\, S^{n/2},
\eeq
where $S$ (the best Sobolev constant) is defined by
\beq
S=\inf\left\{\int_{\R^n}|\D u|^2dx\ :\ u\in H^1(\R^n),\
  \int_{\R^n}|u|^{2n\over n-2}dx=1\right\}.
\eeq
\end{teo}
\no Let us point out that Theorem \ref{Tn3.4} gives a new result also
when $w\equiv0$ in $\Omega$. 
In fact, in this case the functional $\cF$ is even but well known
results (see \cite{BN83,CSS86,S84}) guarantee only the existence of a
finite number 
of solutions (because some compactness conditions hold only at
suitable levels of $\cF$).
On the contrary our method,  combined with some estimates as in
\cite{BN83} and in \cite{CSS86}, allows us to construct infinitely many
solutions with many nodal regions and arbitrarily large energy level.


\vspace{2mm}

{\small {\bf Funding}. 
The authors have been supported by the ``Gruppo
Nazionale per l'Analisi Matematica, la Probabilit\`a e le loro
Applicazioni (GNAMPA)'' of the {\em Istituto Nazionale di Alta Matematica
(INdAM)}.

R.M. acknowledges also the MIUR Excellence Department
Project awarded to the Department of Mathematics, University of Rome
Tor Vergata, CUP E83C18000100006. 
 }



\end{document}